# On the Normality, Regularity and Chain-completeness of Partially Ordered Banach Spaces and Applications

Jinlu Li
Department of Mathematics
Shawnee State University
Portsmouth, Ohio 45662
USA

**Abstract**

In this paper, we study the connections between the normality, regularity, full regularity, and chain-complete property in partially ordered Banach spaces. Then, by applying these properties, we prove some fixed point theorems on partially ordered Banach spaces. As applications of these fixed point theorems, we prove the existence of solutions of some integral equations, such as Hammerstein integral equations, in Banach spaces.

**1. Introduction**

The normality, regularity, full regularity and chain-complete property in partially ordered Banach spaces play important roles in analysis on Banach spaces. They have been widely applied to fixed point theory on partially ordered sets, on which topological structures are either equipped or not (see [1], [9], [11], [15], [19]). The connections between these properties have been studied by some authors (see [10]). For example, from Lemma 2.3.1 in Guo [10], we obtain the following result immediately:

**Lemma**. *Every closed order-interval in a regular partially ordered Banach space is chain-complete*.

In section 3, these properties will be more deeply studied. It is well known that in the traditional fixed point theory the underlying spaces are topological spaces and that considered mappings must satisfy a certain type of continuity with respect to the topologies on the underlying spaces. In contrast with the traditional fixed point theory, for fixed point theorems on lattices or posets (such as the above listed fixed point theorem), the underlying spaces only hold some ordering structures (lattice or partial order), in which, there is neither algebraic structure nor topological structure. The considered mappings are only required to satisfy some ordering increasing conditions. Therefore, fixed point theorems on lattices or posets provide powerful tools to study ordering problems without any topological structure, such as ordered optimizations problems and generalized Nash Equilibrium problems for strategic games with incomplete preferences (see [1-2], [5], [8], [11-12], [15], [21], [23]).

As a topological space or a Banach space is equipped with a partial order, the topological space or the Banach space becomes a partially ordered topological space or partially ordered Banach space. Then, in the study about Banach spaces, in addition to the original algebra structure and the topological structure, the ordering structure will provide an additional tool. The ordering relation on Banach spaces will play important roles to remedy the defect of the non-continuity of the considered mappings and the non-compactness of the underlying spaces. More precisely speaking, when one tries to solve a problem in Banach spaces, it is possible that the continuity condition of the considered mappings and the compactness condition of the underlying spaces may be replaced by some partial ordering properties. In this respect, fixed point theorems have shown the power and effectiveness in solving ordinary differential equations, integral equations, nonlinear fractional evolution equations, and ordered variatinal inequalities (See [3], [6-7], [10], [13-14], [16-18], [22]).

In fixed point theory on partially ordered Banach spaces, some conditions or properties of the given partial order must be satisfied to guarantee the existence of a fixed point for a considered mapping. These properties include normality, regularity, full regularity and chain-complete, etc. The connections between the normal and regular properties on partially ordered Banach spaces have been studied by various authors. For example, in [10], Guo provided the characteristics of the normality, regularity, and full regularity and their inclusion properties on partially ordered Banach spaces. By using these properties, Guo proved several fixed point theorems. With these fixed point theorems, numerous solution existence theorems for ordinary differential equations and integral equations were proven.

In this paper, we will more deeply study the connections between the normality, regularity, full regularity and chain-complete properties on partially ordered Banach spaces. We independently show that regularity implies chain-complete. By applying Theorem 3.1 in [15], we provide some fixed point theorems on both partially ordered Banach spaces and regular partially ordered Banach spaces. Then using these fixed point theorems, we prove the solvability of some integral equations in which the involved functions are not required to satisfy any type of continuity.

## 2. Preliminaries

### 2.1 Natural topology on partially ordered topological space

Let $(X, \succcurlyeq)$ be a poset. For any $u, w \in X$, the following $\succcurlyeq$-intervals are defined:

$$[u) = \{x \in X: x \succcurlyeq u\}, (w] = \{x \in X: x \preccurlyeq w\} \text{ and } [u, w] = [u) \cap (w] = \{x \in X: u \preccurlyeq x \preccurlyeq w\}.$$

If every chain of a subset $D$ of $(X, \succcurlyeq)$ has the smallest $\succcurlyeq$-upper bound in $D$, then $D$ is said to be chain-complete. If every chain of a subset $D$ of $(X, \succcurlyeq)$ has both the smallest $\succcurlyeq$-upper bound and the greatest $\succcurlyeq$-lower bound, then $D$ is said to be bi-chain-complete.

A poset $(X, \succcurlyeq)$ is said to have *the chain-complete property* (or it is simply said to be *chain-complete*) if every $\succcurlyeq$-interval $(w]$ of $(X, \succcurlyeq)$ has is chain-complete. It is said to have *the bi-chain-complete property* (or it is simply said to be *bi-chain-complete*) if every $\succcurlyeq$-interval $[u, w]$ is bi-chain-complete.

Let $(X, \succcurlyeq)$ be a poset equipped with a topology $\tau$ ($(X, \tau)$ is a topological space). The topology $\tau$ is called a natural topology on $(X, \succcurlyeq)$ with respect to the partial order $\succcurlyeq$, whenever, for every $u \in P$, the intervals $[u)$ and $(u]$ both are $\tau$-closed.

In this paper, a poset $(X, \succcurlyeq)$ equipped with a natural topology $\tau$ with respect to $\succcurlyeq$ on $X$ is called a partially ordered topological space; and it is denoted by $(X, \tau, \succcurlyeq)$.

Let $(X, \succcurlyeq)$ be a poset equipped with a metric $d$ (it is also a topological space with the topology on $X$ induced by $d$). If the $d$-topology on $X$ is natural with respect to the partial order $\succcurlyeq$, then $(X, d, \succcurlyeq)$ is called a partially ordered metric space.

Let $X$ be a real vector space endowed with a partial order $\succcurlyeq$, a poset. If the following (order-linearity) properties hold:

1. $x \succcurlyeq y$ implies $x + z \succcurlyeq y + z$, for all $x, y, z \in X$.
2. $x \succcurlyeq y$ implies $\alpha x \succcurlyeq \alpha y$, for all $x, y \in X$ and $\alpha \geq 0$,
3. there are distinct points $x, y \in X$ satisfying $x \succcurlyeq y$, (1)

then $X$ is called a partially ordered vector space. The positive cone of a partially ordered vector space $(X, \succcurlyeq)$ is denoted by $X^+$, which is a convex cone of $X$, and

$$X^+ = [0) = \{x \in X: x \succcurlyeq 0\}.$$

On the other hand, for an arbitrary, nonempty convex cone $K$ of a real vector space $X$, we can define an ordering relation $\succcurlyeq$ on $X$ by

$$x \succcurlyeq y \text{ if and only if } x - y \in K, \text{ for all } x, y \in X.$$

Then $\succcurlyeq$ is a partial order on $X$ induced by $K$ and $(X, \succcurlyeq)$ is a partially ordered vector space, in which the positive cone $X^+ = K$. In this sense, the partial order $\succcurlyeq$ in a partially ordered vector space $(X, \succcurlyeq)$ can be considered as the ordering relation induced by its positive cone $X^+$. A partially ordered topological vector space is both a partially ordered topological space and a partially ordered vector space.

Furthermore, if $X$ is a topological vector space equipped with a partial order $\succcurlyeq$ induced by a cone $K$, then $(X, \succcurlyeq)$ is a partially ordered topological vector space (the topology is natural with respect to $\succcurlyeq$) if and only if the cone $K$ is convex and closed.

Particularly, a Banach space equipped with a partial order is called a partially ordered Banach space if the norm-topology of this space is natural with respect to the given partial order and it satisfies the order-linearity properties in (1).

It is well known that the norm topology of a Banach lattice is always natural with respect to the given lattice order. Hence any Banach lattice can be considered as a special case of partially ordered Banach spaces.

**Remarks**: In this paper, when we refer to the positive cone $X^+$ of a partially ordered topological vector space $(X, \tau, \succcurlyeq)$ it is always $\tau$-closed.

We recall some results from [12] about chain-complete subsets in partially ordered topological spaces.

**Lemma 2.1.** *Let $(X, \tau, \succcurlyeq)$ be a partially ordered Hausdorff topological space. Then every nonempty compact subset of $X$ is chain complete.*

Notice that a subset of a Banach space is (norm) strongly closed if and only if it is weakly closed. It implies that the norm topology in a partially ordered Banach space $(X, \succcurlyeq)$ is natural with respect to $\succcurlyeq$ if and only if the weak topology is natural with respect to $\succcurlyeq$. From Lemma 2.1, it follows:

**Lemma 2.2.** *Every non-empty norm-bounded closed and convex subset of a partially ordered reflexive Banach space is chain-complete.*

## 2.2 The connections between topological limits and order-superiors of sequences in partially ordered topological spaces

Let $\{x_\alpha\}_{\alpha \in \Gamma}$ be a chain in a partially ordered topological space $(X, \tau, \succcurlyeq)$, where $(\Gamma, \succcurlyeq^\Gamma)$ is the index set of this chain, which is a linear ordered set. A point $x \in X$ is called a limit of the chain $\{x_\alpha\}$ if for any given $\tau$-neighborhood $\Delta$ of $x$, there is an index $\beta$ such that

$$\{x_\alpha : \alpha \succcurlyeq^\Gamma \beta\} \subseteq \Delta. \tag{2}$$

It is clear to see that the definition of limit defined by (2) is a natural extension of the sequential limit in analysis.

**Lemma 2.3.** *Let $\{x_n\}$ be a sequence (a chain) in a partially ordered Hausdoff topological space $(X, \tau, \succcurlyeq)$. Suppose that $x$ is the limit of $\{x_n\}$. We have*

  (a) *if $\{x_n\}$ is $\succcurlyeq$-increasing, then $\vee\{x_n\} = x$;*
  (b) *if $\{x_n\}$ is $\succcurlyeq$-decreasing, then $\wedge\{x_n\} = x$.*

*Proof.* We only prove part (a). Part (b) can be similarly proved. For any $m \in N$, where $N$ denotes the set of positive integers, since $[x_m)$ is $\tau$-closed, $\{x_n\}$ is $\succcurlyeq$-increasing, and $x$ is the limit of $\{x_n\}$, it implies $x \in [x_m)$; and therefore,

$$x_m \preccurlyeq x, \text{ for all } m \in N.$$

It follows that $x$ is an upper bound of $\{x_n\}$. To show that $x$ is the least $\succcurlyeq$-upper bound of $\{x_n\}$, let $y$ be any given $\succcurlyeq$-upper bound of $\{x_n\}$. Since $\{x_n\} \subseteq (y]$, then from the $\tau$-closeness of $(y]$ and $x$ being a limit of $\{x_n\}$, it implies $x \in (y]$. Hence $x \preccurlyeq y$. □

The following lemma and Example 2.7 show that the converse of Lemma 2.3 does not hold.

**Lemma 2.4**. *There are some partially ordered non-reflexive Banach spaces $(X, \succcurlyeq)$, in which there is a sequence $\{x_n\}$ such that*

$$\vee\{x_n: 1 \leq n < \infty\} \text{ exists} \quad \text{and} \quad \{x_n\} \text{ is not a Cauchy sequence.}$$

*That is,*

$$\vee\{x_n: 1 \leq n < \infty\} = v, \text{ for some } v \in X, \text{ but } \lim x_n \text{ does not exist.}$$

Proof. We will construct a counter example. Let $C[0, 2]$ be the Banach space of all continuous real valued functions on $[0, 2]$ with the absolute maximum norm ($\|x\| = \max_{0 \leq t \leq 2} |x(t)|$, for any $x \in C[0, 2]$). It is a non-reflexive Banach space. Let $\succcurlyeq$ be the component-wise partial ordering on $C[0, 2]$. That is, for $x, y \in C[0, 2]$,

$$y \succcurlyeq x \text{ if, and only if } y(t) \geq x(t), \text{ for all } t \in [0, 2].$$

It can be shown that the topology induced by the norm $\|\cdot\|$ is natural with respect to the partial order $\succcurlyeq$ on $C[0, 2]$. Hence $(C[0, 2], \|\cdot\|, \succcurlyeq)$ is a partially ordered non-reflexive Banach space. For any $n$, define $x_n \in C(0, 2)$ by

$$x_n(t) = \begin{cases} 0, & \text{if } t = 0, \\ nt, & \text{if } t \in \left(0, \dfrac{1}{n}\right), \\ 1, & \text{if } t \in \left[\dfrac{1}{n}, 2\right]. \end{cases}$$

Define a non-continuous function $u$ on $[0, 2]$ as follows:

$$u(t) = \begin{cases} 0, & \text{if } t = 0, \\ 1, & \text{if } t \in (0, 2]. \end{cases}$$

Since for all $n$, $\|x_n\| = 1$, we have $x_n(t) \uparrow u(t)$, for all $t \in [0, 2]$, as $n \to \infty$, but $\{x_n\}$ is not a Cauchy sequence in $(C[0, 2], \|\cdot\|, \succcurlyeq)$. On the other hand, it can be seen that, for all $n$, $x_n \preccurlyeq x_{n+1}$; therefore $\{x_n\}$ is an $\succcurlyeq$-increasing sequential chain in $(C[0, 2], \|\cdot\|, \succcurlyeq)$. Let $v(t) = 1$, for all $t \in [0, 2]$. Then $v \in C[0, 2]$ and one can check that

$$\vee\{x_n\} = v,$$

But $\{x_n\}$ is not a Cauchy sequence in $(C[0, 2], \|\cdot\|, \succcurlyeq)$. □

**Example 2.7**. Let $\succcurlyeq$ be the component-wise partial ordering on $l_\infty$, that is, for $x = \{s_m\}, y = \{t_m\} \in l_\infty$,

$$y \succcurlyeq x \text{ if, and only if } t_m \geq s_m, \text{ for } m = 1, 2, \ldots .$$

The topology induced by the supremum norm $\|\cdot\|$ is natural with respect to the partial order $\succcurlyeq$ on $l_\infty$. Hence, $(l_\infty, \|\cdot\|, \succcurlyeq)$ is a partially ordered non-reflexive Banach space. For every $n$, define $x_n = \{s_{nm}\} \in l_\infty$ as follows

$$s_{nm} = \begin{cases} 1, & m \le n, \\ 0, & m > n. \end{cases}, \text{ for } m = 1, 2, \ldots.$$

Let $w = \{p_m\} \in l_\infty$ with $p_m = 1$, for all $m = 1, 2, \ldots$. Then

$$\vee\{x_n: 1 \le n < \infty\} = w.$$

It is clear to see that $\{x_n\}$ is not a Cauchy sequence in $l_\infty$; therefore,

the limit of $\{x_n\}$ does not exist as $n \to \infty$, (the norm topology on $l_\infty$). □

## 2.3 Partially ordered non-reflexive Banach spaces may not have the chain-complete property

The Lemma 2.2 in the last subsection states that every partially ordered reflexive Banach space has the chain-complete property. That is, every non-empty norm-bounded closed and convex subset of a partially ordered reflexive Banach space is chain-complete. In this subsection, we show that some partially ordered non-reflexive Banach spaces do not have the chain-complete property (See Lemma 2.8 below). Furthermore, we show that in some partially ordered non-reflexive Banach spaces, there are some norm-bounded and closed subsets that are not even inductive (See Lemma 2.9 below).

**Lemma 2.8.** *A non-empty order-bounded, closed and convex subset of a partially ordered non-reflexive Banach space may be both inductive and not chain-complete.*

*Proof.* We provide an example to prove this lemma. Let $(C[0, 2], \|\cdot\|, \succcurlyeq)$ be the partially ordered non-reflexive Banach space defined in Lemma 2.4. Let $B(0, 1) = \{x \in C[0, 2]: \|x\| \le 1\}$. Let $v$, and $u$ be defined as constant functions with values 1 and $-1$, respectively. For every $n$, define

$$x_n(t) = \begin{cases} 0, & \text{if } t \in [0, 1], \\ nt - n, & \text{if } t \in (1, 1 + \frac{1}{n}], \\ 1, & \text{if } t \in (1 + \frac{1}{n}, 2]. \end{cases}$$

Then, for all $n$, $x_n \in C[0, 2]$ and $\|x_n\| = 1$. It is clear to see that $x_n(t) \uparrow w(t)$, as $n \to \infty$, for all $t \in [0, 2]$, where

$$w(t) = \begin{cases} 0, & \text{if } t \in [0, 1], \\ 1, & \text{if } t \in (1, 2]. \end{cases}$$

Notice that $w \notin C[0, 2]$ and $\{x_n\}$ is not a Cauchy sequence in $C[0, 2]$ (It does not converge to $w$ in $C[0, 2]$). For all $n$, $x_n \preccurlyeq x_{n+1}$ and $\{x_n\}$ is an $\succcurlyeq$-increasing chain in $(C[0, 2], \|\cdot\|, \succcurlyeq)$. One can check that

$$\vee\{x_n: 1 \leq n < \infty\} \text{ does not exist.}$$

It implies that the $\succcurlyeq$-bounded closed and convex subset $(B(0, 1), \succcurlyeq) \subseteq (C[0, 2], \|\cdot\|, \succcurlyeq)$ is not chain-complete.

Since $v$ and $u$ are the $\succcurlyeq$-largest and $\succcurlyeq$-least points in $B(0, 1)$, respectively. So $(B(0, 1), \succcurlyeq)$ is a bi-inductive, $\succcurlyeq$-bounded closed and convex subset of $(C[0, 2], \|\cdot\|, \succcurlyeq)$. □

Furthermore, we have the following result.

**Lemma 2.9**. *A non-empty norm-bounded closed and convex subset of a partially ordered non-reflexive Banach space doesn't need to be inductive.*

*Proof.* We will construct a counter example. Let $C^1[-1, 1]$ denote the space of continuously differentiable functions on $[-1, 1]$. A norm on $C^1[-1, 1]$ is defined as:

$$\|x\| = \max_{-1 \leq t \leq 1} |x(t)| + \max_{-1 \leq t \leq 1} |x'(t)|, \text{ for } x \in C^1[-1, 1]. \tag{3}$$

It is well known that $C^1[-1, 1]$ is a non-reflexive Banach space with the above norm. Define a partial order $\succcurlyeq$ on $C^1[-1, 1]$: for $x, y \in C^1[-1, 1]$, as

$$y \succcurlyeq x \text{ if and only if } y(t) \geq x(t) \text{ and } y'(t) \geq x'(t), \text{ for all } t \in [0, 2]. \tag{4}$$

It can be shown that the topology of the norm $\|\cdot\|$ defined in (3) is natural with respect to the partial order $\succcurlyeq$ on $C^1[-1, 1]$. Hence $(C^1[-1, 1], \|\cdot\|, \succcurlyeq)$ is a partially ordered non-reflexive Banach space.

Let $B(0, 2) = \{x \in C^1[-1, 1]: \|x\| \leq 2\}$, which is a non-empty $\|\cdot\|$-bounded closed and convex subset of $C^1[-1, 1]$. For any given $\lambda \in (0, 1)$, let

$$y(t) = \begin{cases} t + 1, & \text{if } t \in [-1, -\lambda], \\ -\dfrac{1}{2\lambda}t^2 + 1 - \dfrac{\lambda}{2}, & \text{if } t \in (-\lambda, 0], \\ 1 - \dfrac{\lambda}{2}, & \text{if } t \in (0, 1]. \end{cases} \tag{5}$$

It is clear to see that $y \in C^1[-1, 1]$ with $\|y\| = 2 - \lambda/2$. Take a strictly decreasing sequence $\{\lambda_n\} \subseteq (0, 1)$ satisfying

$$\lambda_{n+1} < \frac{\lambda_n}{2}, \text{ for } n = 1, 2, 3, \ldots. \tag{6}$$

For every $n$, we define an element $y_n \in C^1[-1, 1]$ similarly as the definition of the function $y$ in (5), where $\lambda_n$ is substituted for $\lambda$. Let

$$v(t) = \begin{cases} t+1, & \text{if } t \in [-1, 0), \\ 1, & \text{if } t \in [0, 1]. \end{cases}$$

Then

1. $v \notin C^1[-1, 1]$;
2. $y_n \in C^1[-1, 1]$ with $\|y_n\| = 2 - \dfrac{\lambda_n}{2}$, for $n = 1, 2, 3, \ldots$;
3. for every $t \in [-1, 1]$, $y_n(t) \uparrow v(t)$, as $n \to \infty$;
4. for every $n$, $y_n \preccurlyeq y_{n+1}$. That is $\{y_n\}_{1 \leq n < \infty}$ is an $\succcurlyeq$-increasing chain in $(B(0, 2), \succcurlyeq)$;
5. $\{y_n\}_{0 \leq n < \infty}$ is not a Cauchy sequence in $B(0, 2)$;
6. $\{y_n : 1 \leq n < \infty\}$ does not have an $\succcurlyeq$-upper bound in $(B(0, 2), \succcurlyeq)$.

It implies that the non-empty $\|\cdot\|$-bounded closed and convex subset $(B(0, 2), \succcurlyeq)$ of the non-reflexive Banach space $C^1[-1, 1]$ is not inductive.

We only prove the last two parts.

Proof of part 5. For $n = 1, 2, 3, \ldots$, from (6), we have

$$y'_{n+1}(-\lambda_{n+1}) - y'_n(-\lambda_{n+1}) = 1 - \frac{\lambda_{n+1}}{\lambda_n} > \frac{1}{2}.$$

It implies

$$\|y_{n+1} - y_n\| > \frac{1}{2}, \text{ for } n = 1, 2, 3, \ldots.$$

Hence $\{y_n\}_{0 \leq n < \infty}$ is not a Cauchy sequence in $B(0, 2) \subseteq C^1[-1, 1]$.

Proof of part 6. For $n = 1, 2, 3, \ldots$, from definition (3), we have

$$\max_{-1 \leq t \leq 1} |y_n(t)| = \max_{-1 \leq t \leq 1} y_n(t) = 1 - \frac{\lambda_n}{2},$$

and

$$\max_{-1 \leq t \leq 1} |y'_n(t)| = \max_{-1 \leq t \leq 1} y'_n(t) = 1.$$

More precisely, we have

$$y_n(t) = 1 - \frac{\lambda_n}{2}, \text{ for } 0 \leq t \leq 1, \tag{7}$$

and

$$y'_n(t) = 1, \text{ for } -1 \leq t \leq -\lambda_n. \tag{8}$$

Assume, by the way of contradiction, that $\{y_n: 1 \leq n < \infty\}$ has an $\succcurlyeq$-upper bound $w$ in the subposet $(B(0, 2), \succcurlyeq)$. Then

$$w \succcurlyeq y_n, \text{ for } 1 \leq n < \infty.$$

From the definition of the partial order $\succcurlyeq$ in (4), for $1 \leq n < \infty$, we have

$$w(t) \geq y_n(t) \text{ and } w'(t) \geq y'_n(t), \text{ for all } t \in [0, 2].$$

From (7) and (8), we obtain

$$w(t) \geq 1, \text{ for } 0 \leq t \leq 1, \tag{9}$$

and

$$w'(t) \geq 1, \text{ for } -1 \leq t < 0. \tag{10}$$

On the other hand, from $w \in B(0, 2)$ and from (3), (9) and (10), we have

$$2 \geq \|w\| = \max_{-1 \leq t \leq 1}|w(t)| + \max_{-1 \leq t \leq 1}|w'(t)| = \max_{0 \leq t \leq 1} w(t) + \max_{-1 \leq t \leq -\lambda_n} w'(t) \geq 2.$$

It implies

$$w(t) = 1, \text{ for } 0 \leq t \leq 1, \tag{11}$$

and

$$w'(t) = 1, \text{ for } -1 \leq t < 0. \tag{12}$$

(12) and (11) imply that $w$ is linear on $[-1, 0)$ with slope 1 and passing through point $(-1, 0)$, and $w(t) = 1$, for $0 \leq t \leq 1$. So $w$ must coincide with $v$ defined in (3). It is a contradiction to the fact that $v \notin C^1[-1, 1]$.

From this part 6, we conclude that a non-empty norm-bounded closed and convex subset of a partially ordered non-reflexive Banach space may not be inductive.

As a consequence of Lemma 2.9, we obtain similar statement for chain-completeness

**Lemma 2.10**. *A non-empty norm-bounded closed and convex subset of a partially ordered non-reflexive Banach space doesn't need to be chain-complete*.

### 2.4 Normal and regular partially ordered Banach spaces

In this subsection, we recall some concepts of normal and regular partially ordered Banach spaces. For details, the readers are referred to Guo [10].

Let $(X, \succcurlyeq, \|\cdot\|)$ be a partially ordered Banach space. If there is a constant $\lambda > 0$ such that

$$0 \preccurlyeq x \preccurlyeq y \text{ implies that } \|x\| \leq \lambda\|y\|, \tag{13}$$

then ⩾ is said to be normal and $(X, \geqslant, \|\cdot\|)$ is called a normal partially ordered Banach space. The minimum value of $\lambda$ satisfying (13) is called the normal constant of ⩾. It is clear that this minimum value $\lambda$ exists and is greater than or equal to 1.

If every ⩾-upper bounded ⩾-increasing sequence $\{x_n\}$ of $X$ is an $\|\cdot\|$-convergent sequence, that is,

$$x_1 \leqslant x_2 \leqslant \ldots \leqslant y, \text{ for some } y \in X \implies \text{there is } x \in X \text{ such that } \|x_n - x\| \to 0 \text{ as } n \to \infty,$$

then ⩾ is said to be regular and $(X, \geqslant, \|\cdot\|)$ is called a regular partially ordered Banach space.

If every $\|\cdot\|$-bounded ⩾-increasing sequence $\{x_n\}$ of $X$ is an $\|\cdot\|$-convergent sequence, that is,

$$x_1 \leqslant x_2 \leqslant \ldots, \text{ and } \|x_n\| \leq M < \infty, \text{ for some } M > 0 \implies \|x_n - x\| \to 0 \text{ as } n \to \infty, \text{ for some } x \in X,$$

then ⩾ is said to be fully regular and $(X, \geqslant, \|\cdot\|)$ is called a fully regular partially ordered Banach space.

**Theorem 2.2.2** in [10]. *Let $(X, \|\cdot\|, \geqslant)$ be a partially ordered Banach space. Then*

$$\geqslant \text{ is fully regular} \implies \geqslant \text{ is regular} \implies \geqslant \text{ is normal}.$$

**Theorem 2.1.1** in [10]. *Let $(X, \|\cdot\|, \geqslant)$ be a partially ordered Banach space. Then the following statements are equivalent*:

(i)   ⩾ *is normal*;
(ii)  *the norm $\|\cdot\|$ has an equivalent norm $\|\cdot\|_1$ such that $0 \leqslant x \leqslant y$ implies $\|x\|_1 \leq \|y\|_1$*;
(iii) *every ⩾-interval $[x, y] = \{z \in X : x \leqslant z \leqslant y\}$ is $\|\cdot\|$-bounded*;

**Question**: In a regular partially ordered Banach space $(X, \|\cdot\|, \geqslant)$, suppose that $\|\cdot\|$ satisfies

$$0 \leqslant x \leqslant y \text{ implies that } \|x\| \leq \|y\|, \text{ for } x, y \in X^+.$$

Does $0 \leqslant x < y$ imply $\|x\| < \|y\|$, for $x, y \in X^+$ ? The answer for the above question is no, which is shown by the following lemma.

**Lemma 2.11**. *There are some regular partially ordered Banach spaces $(X, \|\cdot\|, \geqslant)$, in which there are some points $x, y \in X^+$ such that*

$$0 \leqslant x < y \text{ and } \|x\| = \|y\|.$$

*Proof.* In $R^2$, define $\|\cdot\|$ to be the sum of absolute values of the coordinates. Let ⩾ be the partial order on $R^2$ induced by the following cone:

$$K = \{(s, t) \in R^2 : |s| \leq t\}$$

That is, for any $(x_1, y_1), (x_2, y_2) \in R^2$,

$$(x_1, y_1) \leqslant (x_2, y_2) \text{ if and only if } |x_2 - x_1| \leq y_2 - y_1.$$

It is enclosed by the two bisectors of QI and QII in $R^2$. Since $K$ is a closed convex cone in $(R^2, \|\cdot\|)$, the topology of the norm $\|\cdot\|$ is natural with respect to the partial order ⩾ induced by $K$.

Let $(x_n, y_n)$, $n = 1, 2, \ldots$, be an $\succcurlyeq$-increasing sequence such that the norms are bounded by a positive number $M$:

$$|x_m - x_n| \leq y_m - y_n, \text{ for all } n \leq m, \tag{14}$$

and

$$|x_n| + |y_n| \leq M, \text{ for all } n. \tag{15}$$

From (14) and (15), it implies that $\{y_n\}$ is an increasing sequence of nonnegative numbers with an upper bound $M$. Hence it is convergent, and so it is a Cauchy sequence. Then, from (14), $\{x_n\}$ is also a Cauchy sequence. It implies that $\{(x_n, y_n)\}$ is a Cauchy sequence in $(R^2, \|\cdot\|)$ and $(R^2, \succcurlyeq, \|\cdot\|)$ is a fully regular partially ordered Banach space. We have

$$0 \preccurlyeq \left(-\frac{1}{2}, \frac{1}{2}\right) \prec (0, 1) \quad \text{and} \quad \left\|\left(-\frac{1}{2}, \frac{1}{2}\right)\right\| = \|(0, 1)\| = 1.$$

All points on the line segment connecting points $\left(-\frac{1}{2}, \frac{1}{2}\right)$ and $(0, 1)$ form an $\succcurlyeq$-increasing chain, in which all points have the same norm 1. Also

$$0 \preccurlyeq \left(\frac{1}{2}, \frac{1}{2}\right) \prec (0, 1) \quad \text{and} \quad \left\|\left(\frac{1}{2}, \frac{1}{2}\right)\right\| = \|(0, 1)\| = 1.$$

All points on the line segment connecting points $\left(\frac{1}{2}, \frac{1}{2}\right)$ and $(0, 1)$ form an $\succcurlyeq$-increasing chain, in which all points have the same norm 1. □

**Lemma 2.12**. *The normality of a partially ordered Banach space does not guarantee the chain-complete property.*

*Proof.* Let $(C[0, 2], \|\cdot\|, \succcurlyeq)$ be the partially ordered non-reflexive Banach space defined in Lemma 2.4. Let $u, v$ be the constant functions with values $-1$ and $1$, respectively. It is straightforward to check that $(C[0, 2], \succcurlyeq, \|\cdot\|)$ is a normal partially ordered Banach space. (Note that from the proof of Lemma 2.4, we see that $(C[0, 2], \succcurlyeq, \|\cdot\|)$ is not regular). From the proof of Lemma 2.8, we get that the $\succcurlyeq$-interval $[u, v]$ is not chain-complete. It concludes that $(C[0, 2], \succcurlyeq, \|\cdot\|)$ does not have the chain-complete property. □

In Theorem 3.1 of the later sections, we will show that every nonempty closed inductive subset of a regular partially ordered Banach space is chain complete. As a consequence, it implies that the regularity of partially ordered Banach spaces guaranties that every closed order-interval is chain-complete.

## 3. Chain-complete property and fixed point on regular partially ordered Banach spaces

### 3.1 Regular partially ordered Banach spaces have chain-complete property

In the previous section, we showed that some closed order-intervals (which are special inductive subsets) of some normal partially ordered Banach spaces are not chain complete. In this section, we prove that every regular partially ordered Banach space has the chain-complete property, which is showed by proving the fact that every closed inductive subset of any given regular partially ordered Banach space are chain-complete. Then, by applying this property, we prove several fixed point theorems of set-valued mappings on regular partially ordered Banach spaces.

The following theorem can be proved by using Lemma 2.3.1 in Guo [10]. Here we give a direct and simpler proof.

**Theorem 3.1**. *Every nonempty closed inductive subset of a regular partially ordered Banach space is chain-complete*.

*Proof.* Let $(X, \|\cdot\|, \succcurlyeq)$ be a regular partially ordered Banach space and let $D$ be an arbitrary nonempty closed inductive subset of $X$. Let $C = \{x_\alpha\}$ be an arbitrary ($\succcurlyeq$-increasing with respect to the index ordering) chain in $D$. Since $D$ is inductive, then $\{x_\alpha\}$ has an $\succcurlyeq$-upper bound $w \in D$. That is,

$$x_\alpha \preccurlyeq w, \text{ for all index } \alpha. \tag{16}$$

Take an arbitrary fixed index $\beta$. Define $C_\beta = \{x_\alpha\}_{\alpha \in [\beta)}$, that is also an ($\succcurlyeq$-increasing) chain in $D$ and $w \in D$ is also an $\succcurlyeq$-upper bound of $C_\beta$. Let $\overline{C}_\beta$ be the closure of $C_\beta = \{x_\alpha\}_{\alpha \in [\beta)}$. Next we show that $\overline{C}_\beta$ is also a chain in $D$.

Since $(X, \|\cdot\|, \succcurlyeq)$ is regular, it is normal. From the properties of normal partially ordered Banach spaces, we may assume that

$$0 \preccurlyeq x \preccurlyeq y \text{ implies that } \|x\| \le \|y\|, \text{ for } x, y \in X. \tag{17}$$

For any different $x_\gamma, x_\lambda \in \overline{C}_\beta$, let $\delta = \|x_\gamma - x_\lambda\| > 0$. Then we can take two sequences $\{u_n\}$ and $\{v_m\}$ in $C_\beta$ such that

$$u_n \to x_\gamma, \text{ as } n \to \infty, \text{ and } v_m \to x_\lambda, \text{ as } m \to \infty, \tag{18}$$

and satisfying

$$\|x_\gamma - u_n\| < 0.1\delta \text{ and } \|x_\lambda - v_m\| < 0.1\delta, \text{ for all } m, n = 1, 2, \ldots. \tag{19}$$

Then we show that either $u_n \succcurlyeq v_m$, for all $m, n = 1, 2, \ldots$, or $v_m \succcurlyeq u_n$, for all $m, n = 1, 2, \ldots$. As a matter of fact, assume, by the way of contradiction, that there are numbers $m, n$ and $k$ such that

$$u_k \preccurlyeq v_m \preccurlyeq u_n.$$

Then

$$0 \preccurlyeq v_m - u_k \preccurlyeq u_n - u_k.$$

From (17) we have

$$\|v_m - u_k\| \le \|u_n - u_k\|. \tag{20}$$

Estimating the right side of (20) by using (19) gets

$$\|u_n - u_k\| \leq \|u_n - x_\gamma\| + \|u_k - x_\gamma\| \leq 0.2\delta.$$

The left side of (20) satisfies

$$\|v_m - u_k\| \geq \|x_\lambda - x_\gamma\| - \|x_\gamma - u_k\| - \|x_\lambda - v_m\| \geq 0.8\delta,$$

which is a contradiction. Hence, without loss of generality, we can assume that $u_n \succcurlyeq v_m$, for all $m$, $n = 1, 2, \ldots$ . Then, for any fixed $n$, $\{v_m\} \subseteq (u_n]$. From the closeness of $(u_n]$ and (18), we obtain that $x_\lambda \in (u_n]$. That is, $x_\lambda \preccurlyeq u_n$, for every given $n$. It implies $\{u_n\} \subseteq [x_\lambda)$. From the closeness of $[x_\lambda)$ and (18), we obtain that $x_\gamma \in [x_\lambda)$. It implies $x_\lambda \preccurlyeq x_\gamma$ and $x_\lambda$ and $x_\gamma$ are $\succcurlyeq$-comparable. Hence $\overline{C}_\beta$ is also a chain in $D$.

Any sequence $\{x_n\} \subseteq \overline{C}_\beta$ contains either an $\succcurlyeq$-increasing or an $\succcurlyeq$-decreasing subsequence $\{x_{n(k)}\}$. It is clear that $w$ is an $\succcurlyeq$-upper bound and $x_\beta$ is an $\succcurlyeq$-lower bound of $\{x_{n(k)}\}$. From the regularity of $(X, \|\cdot\|, \succcurlyeq)$, $\{x_{n(k)}\}$ being either $\succcurlyeq$-increasing or $\succcurlyeq$-decreasing, it is always convergent. It implies that $\overline{C}_\beta$ is compact (with respect to $\|\cdot\|$). So $C_\beta$ is relatively compact, and $C_\beta$ is separable. Hence $C_\beta$ has a countable dense subset $\{y_n\} \subseteq C_\beta$.

Next we construct a subsequence of $\{y_n\}$ to show the existence of $\vee C_\beta$. For a fixed prime number $p$, we choose an $\succcurlyeq$-increasing subsequence $\{y_{m(k)}: k = 1, 2, \ldots\} \subseteq \{y_n\}$ or an element $y_{n(j)} \in \{y_n\}$, for some positive integer $n(j)$, as follows: For $k = 1$, if there is a number $n(1) \leq p$, such that

$$y_{n(1)} \succcurlyeq y_n, \text{ for all } n = 1, 2, \ldots,$$

then stop the process and we have

$$y_{n(1)} = \vee\{y_n: n = 1, 2, \ldots\}.$$

Otherwise, we take $m(1) \geq p$, such that

$$y_{m(1)} \succcurlyeq y_n, \text{ for all } n = 1, 2, \ldots m(1) - 1.$$

After $m(1)$ is chosen (Assume that $n(1)$ does not exist). For $k = 2$, if there is a number $n(2)$ with $m(1) \leq n(2) \leq m(1) + p^2$, such that

$$y_{n(2)} \succcurlyeq y_n, \text{ for all } n = m(1), m(1) + 1, m(1) + 2, \ldots,$$

then stop the process and we have

$$y_{n(2)} \succcurlyeq y_n, \text{ for all } n = 1, 2, \ldots,$$

and

$$y_{n(2)} = \vee\{y_n: n = 1, 2, \ldots\}.$$

Otherwise, we take $m(2) \geq m(1) + p^2$, such that

$$y_{m(2)} \succcurlyeq y_n, \text{ for all } n = m(1), m(1) + 1, m(1) + 2, \ldots, m(2) - 1.$$

It satisfies

$$y_{m(2)} \succcurlyeq y_{m(1)}.$$

It implies

$$y_{m(2)} \succcurlyeq y_n, \text{ for all } n = 1, 2, \ldots, m(2) - 1.$$

Continue the process. If there exists an integer $n(j)$, for some $j$, then

$$y_{n(j)} = \vee \{y_n : n = 1, 2, \ldots\}.$$

Otherwise, we will choose an $\succcurlyeq$-increasing subsequence $\{y_{m(k)}\}$ of $\{y_n\}$. From (16), it satisfies

$$y_{m(k)} \preccurlyeq w, \text{ for all } k = 1, 2, \ldots .$$

Since $(X, \succcurlyeq, \|\cdot\|)$ is regular, it implies that $\{y_{m(k)}\}$ is a convergent sequence of $C_\beta$ and

$$y_{m(k)} \to y, \text{ as } k \to \infty, \text{ for some } y \in \overline{C}_\beta \text{ (Since } \overline{C}_\beta \text{ is closed).}$$

From Lemma 2.3, we have

$$y = \vee \{y_{m(k)} : k = 1, 2, \ldots\}. \tag{21}$$

Next we prove

$$y = \vee \{y_n : n = 1, 2, \ldots\}. \tag{22}$$

For any $n = 1, 2, \ldots$, from the selection of $\{y_{m(k)}\}$, there is a number $k$ such that $n < m(k)$. It implies $y_n \preccurlyeq y_{m(k)} \preccurlyeq y$. Hence $y$ is an upper bound of $\{y_n\}$. From $\{y_{m(k)}\} \subseteq \{y_n\}$ and (21), (22) follows immediately.

Since $\{y_n\} \subseteq (y]$ which is closed and $\{y_n\}$ is dense in $C_\beta$, it implies

$$C_\beta \subseteq \overline{C}_\beta \subseteq (y]. \tag{23}$$

So $y$ is an upper bound of $C_\beta$. By $\{y_n\} \subseteq (y]$ again, combining (23) and (22), we have

$$y = \vee C_\beta. \tag{24}$$

By definition of $C$ and $C_\beta$, we immediately obtain

$$y = \vee C. \qquad \square$$

As a consequence of Theorem 3.1, we have

**Corollary 3.2.** *Every nonempty closed bi-inductive subset of a regular partially ordered Banach space $(X, \|\cdot\|, \succcurlyeq)$ is bi-chain complete.*

From this corollary and Theorem 3.1, we obtain Lemma 2.3.1 in Guo [10] as the following corollary.

**Corollary 3.3**. *Let $(X, \|\cdot\|, \succcurlyeq)$ be a regular partially ordered Banach space $(X, \succcurlyeq, \|\cdot\|)$. For any $u$, $v \in X$ with $u \preccurlyeq v$, we have*

(a) $(u]$ *is chain-complete*;
(b) $[u)$ *is re-chain-complete*;
(c) $[u, v]$ *is bi-chain-complete*.

We conclude this subsection with the following property as a consequence of Corollary 3.3.

**Corollary 3.4**. *Every regular partially ordered Banach space has the bi-chain-complete property*.

From Corollary 3.4, we can obtain the following result about reflexive Banach lattices (see Theorem 4.9 and Corollary 4.10 in Aliprantis and Burkinshaw [2]).

**Corollary 3.5**. *Every reflexive Banach lattice is bi-chain-complete*.

*Proof.* Let $(X, \|\cdot\|, \succcurlyeq)$ be a reflexive Banach lattice. From the definition of Banach lattices, for all $x, y \in X$,

$$y \succcurlyeq x \succcurlyeq 0 \text{ implies } \|y\| \geq \|x\|.$$

It implies that $(X, \|\cdot\|, \succcurlyeq)$ is normal. Since $(X, \|\cdot\|, \succcurlyeq)$ is reflexive, it follows that $(X, \|\cdot\|, \succcurlyeq)$ is regular. Then Corollary follows immediately from Corollary 3.4. □

It is worth to note that it has been proved by Theorem 4.9 and Corollary 4.10 in Aliprantis and Burkinshaw that every reflexive Banach lattice has Dedekind chain-complete property (For more details, see page 187 in Aliprantis and Burkinshaw [2]).

## 3.3 Fixed point theorems on regular partially ordered Banach spaces

Let $(X, \succcurlyeq)$, $(U, \succcurlyeq^U)$ be posets and $T: X \to 2^U \setminus \{\emptyset\}$ a set-valued mapping. $T$ is said to be isotone, or to be order-increasing upward, if $x \preccurlyeq y$ in $X$ implies, for any $z \in Tx$, there is a $w \in Ty$ such that $z \preccurlyeq^U w$. $T$ is said to be order-increasing downward, if $x \preccurlyeq y$ in $X$ implies, for any $w \in Ty$, there is a $z \in Tx$ such that $z \preccurlyeq^U w$. If $T$ is both order-increasing upward and order-increasing downward, then $T$ is said to be order-increasing.

In particular, a single-valued mapping $F$ from a poset $(X, \succcurlyeq)$ to a poset $(U, \succcurlyeq^U)$ is said to be order-increasing whenever $x \preccurlyeq y$ implies $F(x) \preccurlyeq^U F(y)$. $F$ is said to be strictly order-increasing whenever $x \prec y$ implies $F(x) \prec^U F(y)$.

Let $(X, \|\cdot\|, \succcurlyeq)$ be a partially ordered Banach space. Let $D$ be a subset of $X$. Let $T: D \to 2^X \setminus \{\emptyset\}$ be a set-valued mapping. For a point $x \in D$, if $x \in Tx$, then $x$ is called a fixed point of $T$. The collection of all fixed points of $F$ is denoted by $\tilde{F}(T)$.

A nonempty subset $A$ of a poset $(X, \succcurlyeq)$ is said to be *universally inductive* in $X$ whenever any given chain $\{x_\alpha\} \subseteq X$ satisfying that every element $x_\beta \in \{x_\alpha\}$ has an $\succcurlyeq$-upper cover in $A$ has an $\succcurlyeq$-upper bound in $A$.

Some useful universally inductive subsets in posets are provided in [11], which are listed in the following two lemmas for easy reference.

**Lemma 3.7** *Every inductive subset $A$ with a finite number of maximal elements in a chain complete poset is universally inductive.*

**Lemma 3.8** *Every nonempty compact subset of a partially ordered Hausdorff topological space is universally inductive.*

For the easy reference, we recall

**Theorem 3.1 in** [15]. *Let $(P, \succcurlyeq)$ be a chain-complete poset and let $F : P \to 2^P \setminus \{\emptyset\}$ be a set-valued mapping satisfying the following three conditions*:

A1. *$F$ is order-increasing upward.*
A2. *$(F(x), \succcurlyeq)$ is universally inductive, for every $x \in P$.*
A3. *there is an element $y_*$ in $P$ and $v_* \in F(y_*)$ with $y_* \preccurlyeq v_*$.*

*Then*

(i) *$(\mathcal{F}(F), \succcurlyeq)$ is a nonempty inductive poset.*

(ii) *$(\mathcal{F}(F) \cap [y_*), \succcurlyeq)$ is a nonempty inductive poset; and $F$ has an $\succcurlyeq$-maximal fixed point $x_*$ with $x_* \succcurlyeq y_*$.*

By Theorem 3.1 and as a consequence of Theorem 3.1 in [15], we have

**Theorem 3.9**. *Let $(X, \|\cdot\|, \succcurlyeq)$ be a regular partially ordered Banach space. Let $D$ be a closed inductive subset of $X$. Let $T: D \to 2^D \setminus \{\emptyset\}$ be a set-valued isotone mapping with universally inductive values. Suppose that there are points $x_0 \in D$, $x_1 \in Tx_0$ satisfying $x_0 \preccurlyeq x_1$. Then*

    (a) *$\mathcal{F}(T)$ is a nonempty inductive subset of $D$;*

    (b) *$\mathcal{F}(T) \cap [x_0)$ is a nonempty inductive subset of $D$.*

**Corollary 3.10**. *Let $(X, \|\cdot\|, \succcurlyeq)$ be a regular partially ordered Banach space. Let $D$ be a closed inductive subset of $X$. Let $F: D \to D$ be an $\succcurlyeq$-increasing single-valued mapping. Suppose that there is $x_0 \in D$ satisfying $x_0 \preccurlyeq Fx_0$. Then*

    (a) *$\mathcal{F}(F)$ is a nonempty inductive subset of $D$;*

    (b) *$\mathcal{F}(F) \cap [x_0)$ is a nonempty inductive subset of $D$.*

### 3.4 The set of fixed points of mappings may not be a sublattice in Hilbert lattices

In the previous section, we mentioned that the set of fixed points $\mathcal{F}(T)$, for set-valued mapping $T$, and $\mathcal{F}(F)$, for single-valued mapping $F$ are always inductive if $T$ and $F$ satisfy the conditions given in Theorem 3.9 and Corollary 3.10, respectively. In this section, we provide some examples to show that, for some mappings defined on lattices, under the same conditions as in Theorem 3.9 and Corollary 3.10, the set of fixed points may not be a sublatice.

**Remark 3.11**. Let $D$ be an $\succcurlyeq$-chain-complete lattice. Let $T: D \to 2^D \setminus \{\emptyset\}$ be an isotone mapping with $\succcurlyeq$-chain-complete sublattice values. Suppose that there are points $x_0 \in D$, $x_1 \in Tx_0$ satisfying

$$x_0 \preccurlyeq x_1.$$

Then $\mathcal{F}(T)$ is a nonempty inductive subset of $D$ which may not be a sublattice of $D$.

This can be shown by the following example. Let $\succcurlyeq$ be the component-wise partial ordering on the 2-Eculidean space $(R^2, \|\cdot\|)$ that is, for $u = (s_1, t_1)$, $v = (s_2, t_2) \in R^2$,

$$v \succcurlyeq u \quad \text{if, and only if} \quad s_2 \geq s_1 \text{ and } t_2 \geq t_1.$$

Then $(R^2, \|\cdot\|, \succcurlyeq)$ is a Hilbert lattice.

Let $D$ be the closed square in $R^2$ with vertices at $(0, 0)$, $(0, 2)$, $(2, 0)$ and $(2, 2)$. Define a set-valued mapping $T: D \to 2^D \setminus \{\emptyset\}$ as follows: for every $(s, t) \in D$,

$$T(s, t) = \begin{cases} \{(s, s)\}, & \text{if } 0 \leq s < 1, 0 \leq t \leq 2, \\ \{(s, s), (s, s-1)\}, & \text{if } 1 \leq s \leq 2, 0 \leq t \leq 2. \end{cases}$$

One can check that $T$ satisfies all conditions in this remark. From Theorem 3.9, we conclude that $\mathcal{F}(T)$ is a nonempty inductive subset of $D$. More precisely, $\mathcal{F}(T)$ is the union of the closed segment with ending points $(0, 0)$ and $(2, 2)$ and the closed segment with ending points $(1, 0)$ and $(2, 1)$. It is clear to see that $\mathcal{F}(T)$ is not a sublattice of $D$. For example, $(1, 1)$, $(1.5, 0.5) \in \mathcal{F}(T)$, but $(1, 1) \vee (1.5, 0.5) = (1.5, 1) \notin \mathcal{F}(T)$ and $(1, 1) \wedge (1.5, 0.5) = (1, 0.5) \notin \mathcal{F}(T)$.

For single-valued mappings, we have the following result:

**Remark 3.12**. Let $D$ be an $\succcurlyeq$-chain-complete sublattice of a reflexive Banach lattice. Let $F: D \to D$ be an increasing mapping. Suppose that there are point $x_0 \in D$ satisfying $x_0 \preccurlyeq Fx_0$. Then $\mathcal{F}(F)$ is a nonempty inductive subset of $D$ which may not be a sublattice of $D$.

We will show this by the following two examples.

**Example 1**. Let $(R^2, \|\cdot\|, \succcurlyeq)$ be the Hilbert lattice defined in the proof of Remark 3.11. Let $A$ denote the closed segment with ending points $(0, 0)$ and $(1, 1)$ and let $B$ denote the closed segment with ending points $(2, 2)$ and $(3, 3)$. Let

$$D = A \cup B \cup \{(1, 2), (2, 1)\}.$$

Then $D$ is an $\succcurlyeq$-chain-complete sublattice of $(R^2, \|\cdot\|, \succcurlyeq)$ and it is not a convex subset of $(R^2, \|\cdot\|)$. We define $F: D \to D$ as below:

$$F(1, 2) = (1, 2), F(2, 1) = (2, 1)$$

and

$$F(s, t) = \begin{cases} (0,0), & \text{if } (s, t) \in A, \\ (3,3), & \text{if } (s, t) \in B. \end{cases}$$

Then $\tilde{F}(T) = \{(0, 0), (1, 2), (2, 1), (3, 3)\}$. It is a nonempty inductive subset of $D$. But $\tilde{F}(T)$ is not a sublattice of $D$. For example, $(1, 2), (2, 1) \in \tilde{F}(T)$, but $(1, 2) \vee (2, 1) = (2, 2) \notin \tilde{F}(T)$ and $(1, 2) \wedge (2, 1) = (1, 1) \notin \tilde{F}(T)$.

**Example 2**. Let $(R^2, \|\cdot\|, \succcurlyeq)$ be the Hilbert lattice defined in the proof of Remark 3.11. Let $D$ denote the closed 4-gon with vertexes at $(0, 0), (1, 2), (2, 1)$ and $(3, 3)$. Then $D$ is a closed and convex subset of $R^2$ and $(D, \succcurlyeq)$ is a chain complete sublattice of $(R^2, \|\cdot\|, \succcurlyeq)$. Let $C$ denote the closed segment with ending point $(1, 2)$ and $(2, 1)$ and $C_1 = C \setminus \{(1, 2), (2, 1)\}$. Let $A$ denote the triangle with vertexes at $(0, 0), (1, 2), (2, 1)$ excluding set $C$ and let $B$ denote the triangle with vertexes at $(1, 2), (2, 1)$ and $(3, 3)$ excluding set $C$. We define $F: D \to D$ as below:

$$F(1, 2) = (1, 2), F(2, 1) = (2, 1)$$

and

$$F(s, t) = \begin{cases} (0,0), & \text{if } (s, t) \in A, \\ (3,3), & \text{if } (s, t) \in B \cup C_1. \end{cases}$$

Then $\tilde{F}(T) = \{(0, 0), (1, 2), (2, 1), (3, 3)\}$. It is a nonempty inductive subset of $D$. But $\tilde{F}(T)$ is not a sublattice of $D$. For example, $(1, 2), (2, 1) \in \tilde{F}(T)$, but $(1, 2) \vee (2, 1) = (2, 2) \notin \tilde{F}(T)$ and $(1, 2) \wedge (2, 1) = (1, 1) \notin \tilde{F}(T)$.

## 4. Hammerstein integral equations in $L_p$-spaces

In this section, we always assume that $(\Sigma, \mu)$ is a $\sigma$-finite measure space. We consider the following nonlinear Hammerstein integral equation

$$x(t) = \int_\Sigma T(t, s) f(s, x(s)) d\mu(s), \tag{25}$$

where the function $T(t,s)$ is called the kernel of this Hammerstein integral equation (25). We study the existence of solutions of equation (25) in the real Banach space $L_p(\Sigma, \mu)$, with norm $\|\cdot\|_p$, for $p > 1$.

Next, we prove some existence theorems of solutions for equation (25) in the space $L_p(\Sigma, \mu)$ by Theorem 3.1 in [15]. First we prove a similar result by Corollary 3.10 in the regular partially ordered Banach space $L_p(\Sigma, \mu)$ with respect to the component-wise ordering, where single-valued mappings are considered as special cases of set-valued mappings with singleton values.

**Theorem 4.1**. *Let T and f be real functions defined on $\Sigma \times \Sigma$ and on $\Sigma \times R$, respectively, satisfying the following conditions*:

(i) *there are some numbers p, q > 1 with* $\frac{1}{p} + \frac{1}{q} = 1$, *such that*

$$\lambda = \int_\Sigma \left( \int_\Sigma |T(s,t)|^q d\mu(s) \right)^{\frac{p}{q}} d\mu(t) < \infty.$$

(ii) $T(t, s) > 0$, *a.e. for* $t, s \in \Sigma$;

(iii) $f(s, \cdot)$ *is increasing, for every fixed* $s \in \Sigma$ *and* $f(s, 0) > 0$, *a.e. for* $s \in \Sigma$;

(iv) *there is* $\gamma > 0$, *such that* $\int_\Sigma |f(s, x(s))|^p d\mu(s) \le \frac{\gamma^p}{\lambda}$, *for all* $x \in L_p(\Sigma, \mu)$ *with* $\|x\|_p \le \gamma$.

*Then equation (25) has a solution* $x^*$ *in* $L_p(\Sigma, \mu)$ *with* $0 < \|x^*\|_p \le \gamma$. *Moreover, the solution set of equation (25) in* $L_p(\Sigma, \mu)$ *is a nonempty inductive subset of* $\{x \in L_p(\Sigma, \mu): \|x\|_p \le \gamma\}$.

*Proof.* For the given positive number $p > 1$ in this theorem, let

$$K = \{x \in L_p(\Sigma, \mu): x(t) \ge 0, \text{ a.e. for } t \in \Sigma\}.$$

It is clear that $K$ is a closed convex cone in $L_p(\Sigma, \mu)$. The partial order $\succcurlyeq$ on $L_p(\Sigma, \mu)$ induced by the cone $K$ is equivalently defined as follows. For any $x, y \in L_p(\Sigma, \mu)$,

$$y \succcurlyeq x \text{ if and only if } y(t) \ge x(t), \text{ a.e. for } t \in \Sigma. \tag{26}$$

It follows that, both, the norm topology and the weak topology, are natural topologies on $L_p(\Sigma, \mu)$ with respect to the partial order $\succcurlyeq$ induced by the cone $K$. Let $\omega$ denote the weak topology on $L_p(\Sigma)$. So both $(L_p(\Sigma, \mu), \|\cdot\|_p, \succcurlyeq)$ and $(L_p(\Sigma, \mu), \omega, \succcurlyeq_K)$ are partially ordered topological spaces. Since $L_p(\Sigma, \mu)$ is a reflexive Banach space, then $(L_p(\Sigma, \mu), \|\cdot\|_p, \succcurlyeq)$ is a partially ordered reflexive Banach space with the norm topology. Let

$$P = \{x \in L_p(\Sigma, \mu): \|x\|_p \le \gamma\}.$$

Then $P$ is the closed ball in $L_p(\Sigma, \mu)$ at center origin and with radius $\gamma$, which is a bounded closed convex subset of $L_p(\Sigma, \mu)$. From Lemma 2.2, $(P, \succcurlyeq)$ is a chain-complete poset. Define a mapping $F$ as

$$(Fx)(t) = \int_\Sigma T(t,s) f(s, x(s)) d\mu(s), \text{ for all } x \in P.$$

From conditions (i) and (iv) of this theorem, for any $x \in P$, we have

$$\int_\Sigma |(Fx)(t)|^p d\mu(t) = \int_\Sigma \left| \int_\Sigma T(t,s) f(s, x(s)) d\mu(s) \right|^p d\mu(t)$$

$$\leq \int_\Sigma \left( \int_\Sigma |T(t,s)|^q d\mu(s) \right)^{\frac{p}{q}} \left( \int_\Sigma |f(s,x(s))|^p d\mu(s) \right) d\mu(t)$$

$$= \int_\Sigma \left( \int_\Sigma |T(t,s)|^q d\mu(s) \right)^{\frac{p}{q}} d\mu(t) \int_\Sigma |f(s,x(s))|^p d\mu(s)$$

$$= \lambda \int_\Sigma |f(s,x(s))|^p d\mu(s)$$

$$\leq \gamma^p.$$

It implies that $F$ is a mapping from $P$ to $P$. For any $x, y \in P$, $y \succcurlyeq x$ is equivalent to $y(t) \geq x(t)$, a.e. $t \in \Sigma$ by (26). From conditions (ii) and (iii) in this theorem, it follows that

$$(Fy)(t) - (Fx)(t)$$

$$= \int_\Sigma T(t,s) f(s, y(s)) d\mu(s) - \int_\Sigma T(t,s) f(s, x(s)) d\mu(s)$$

$$= \int_\Sigma T(t,s) (f(s, y(s)) - f(s, x(s))) d\mu(s)$$

$$\geq 0.$$

Hence $F\colon P \to P$ is an $\succcurlyeq$-increasing single-valued mapping. Hence $F$ satisfies condition A1 in Theorem 3.1 in [15]. Since every singleton is clearly universally inductive, then $F$ satisfies condition A2 in Theorem 3.1 in [15]. To prove that the mapping $F$ satisfies condition A3 in Theorem 3.1 in [15], we chose the constant function $\theta \in L_p(\Sigma, \mu)$ with value 0. From conditions (ii) and (iii) in this theorem, we have

$$(F\theta)(t) = \int_\Sigma T(t,s) f(s, 0) d\mu(s) > 0, \text{ for all } t \in \Sigma.$$

It follows that $F\theta \succcurlyeq \theta$ and $F$ satisfies condition A3 in Theorem 3.1 in [15] with the constant function $\theta$. Then, from Theorem 3.1 in [15], $F$ has a fixed point, say $x^* \in P$, which satisfies $\|x^*\|_p \leq \gamma$. It is a solution of the Hammerstein integral equation (25). From the above inequality, it follows that the constant function $\theta$ is not a solution of the equation (25). It follows that $x^* \neq \theta$. Then from part (a) in Theorem 3.1 in [15], it immediately follows that the solution set of equation (25) in $L_p(\Sigma, \mu)$ is a nonempty inductive subset of $P$. □

# References

[1] Abian S. and Brown A., A theorem on partially ordered sets with applications to fixed point theorem, *Canad. J. Math.*, (1961), 13: 78–83.

[2] Aliprantis, C., and Burkinshaw, O., *Positive Operators*, Springer, The Netherlands, (2006).

[3] Bhaskar B., and Lakshmikantham T., V: Fixed point theorems in partially ordered metric spaces and applications. *Nonlinear Anal*. 65, 1379–1393, (2006).

[4] Carl S. and Heikkilä S., Fixed point theory in ordered sets and applications, in Differential and Integral Equations to Game Theory, Springer, New York, (2010).

[5] Cartwright E. and Wooders M., On purification of equilibrium in Bayesian games and expose Nash equilibrium, *International Journal of Game Theory,* (2009), 38(1): 127–136.

[6] Deimling Klaus, *Nonlinear Functional Analysis*, Springer-Verlag, New York, Heidelberg, Berlin, Tokyo, (1989).

[7] Dhage B., Partially continuous mappings in partially ordered normed linear spaces and applications to functional integral equations, *Tamkang J. Math.*, (2014), 45: 397–426.

[8] Dhage B., Hybrid fixed point theory in partially ordered normed linear spaces and applications to fractional integral equations, *J. Differ. Equ. Appl.*, (2013), 2: 155–184.

[9] Fujimoto, An extension of Tarski's fixed point theorem and its applications to isotonic complementarity problems, *Math. Program.*, (1984), 28: 116–118.

[10] Guo D., Partial order methods in nonlinear analysis, Shandong Academic Press, (1997).

[11] Li J. L., Several extensions of the Abian–Brown fixed point theorem and their applications to extended and generalized Nash equilibria on chain-complete posets, *J. Math. Anal. Appl.*, (2014), 409:1084–1092.

[12] Li J. L., Fixed point theorems on partially ordered topological vector spaces and their applications to equilibrium problems with incomplete preferences, *Fixed Point Theory and Applications*, (2014), 2014/191.

[13] Li J. L., On the existence of solutions of variational inequalities in Banach spaces, *Journal of Mathematical Analysis and Applications*, (2004), 295:115–126.

[14] Li J. L., A lower and upper bounds version of a variational inequality, *Applied Mathematics Letters*, (2000), 13: 47–51.

[15] Li J. L., Inductive properties of fixed point sets of mappings on posets and on partially ordered topological spaces, *Fixed Point Theory and Applications*, (2015) 2015:211 DOI 10.1186/s13663-015-0461-8.

[16] Li J. L., Existence of Continuous Solutions of Nonlinear Hammerstein Integral Equations Proved by Fixed Point Theorems on Posets, *Journal of Nonlinear and Convex Analysis*, Volume **17**, Number **7**, (2016), 1333–1345.

[17] Nieto J. and Lopez R.R., Contractive mapping theorems in partially ordered sets and applications to ordinary differential equations, *Order* (in press).

[18]  Nieto J. and Lopez R.R., Existence and uniqueness of fixed point in partially ordered sets and applications to ordinary differential equations, *Acta. Math. Sinica* (in press).

[19] Tarski A., A lattice-theoretical fixed point theorem and its applications, *Pacific J. Math.*,(1955), 5: 285–309.

[20] Ward L. E. Jr., Partially ordered topological space. *Proc. Am. Math. Soc*. 5(1), 144–161 (1954)

[21] Xie L. S., Li J. L. and Yang W. S., Order-clustered fixed point theorems on chain-complete preordered sets and their applications to extended and generalized Nash equilibria, *Fixed Point Theory and Applications*, (2013), 2013/1/192.

[22] Yang H., Ravi P. A., Hemant K. N. and Yue L., Fixed point theorems in partially ordered Banach spaces with applications to nonlinear fractional evolution equations, *J. Fixed Point Theory Appl.*，(2016), DOI 10.1007/s11784-016-0316-xc.

[23]  Zhang Congjun and Wang Yuehu, Applications of order-theoretic fixed point theorems to discontinuous quasi-equilibrium problems, *Fixed Point Theory and Applications* (2015) 2015:54 DOI 10.1186/s13663-015-0306-5.